\documentclass{amsart}
\usepackage{a4}
\usepackage{amssymb}
\usepackage{verbatim}

\newtheorem{Th}{Theorem}[section]
\newtheorem{Cor}[Th]{Corollary}
\newtheorem{Lem}[Th]{Lemma}
\newtheorem{Prop}[Th]{Proposition}

\newtheorem{lemma-num}{Lemma}
\newtheorem{claim-num}{Claim}
\newtheorem{cor-num}{Corollary}

\newtheorem*{lem}{Lemma}

\def\gl#1{\mathop{\rm GL}(#1)}

\def\aut#1{\mathop{\rm Aut}(#1)}
\def\inn#1{\mathop{\rm Inn}(#1)}
\def\ia#1{\mathop{\rm IA}(#1)}
\def\Fix{\mathop{\rm Fix}}

\def\id{\mathop{\rm id}}
\def\iat{\mathop{IA}}

\def\inv{^{-1}}
\def\str#1{\langle#1\rangle}

\def\avst#1{ \overline{\mathstrut #1} }

\def\f{\varphi}
\def\a{\alpha}

\def\s{\sigma}
\def\vk{\varkappa}

\def\av#1{\overline{#1}}

\def\cB{{\mathcal B}}
\def\cC{{\mathcal C}}

\def\cV{{\mathcal V}}

\def\Z{\mathbf Z}

\def\bmpi{\pi}

\renewcommand{\le}{\leqslant}
\renewcommand{\ge}{\geqslant}

\def\rank{\mathop{\rm rank}\,}

\def\thetast{\theta}

\begin{document}

\title{Free two-step nilpotent groups whose automorphism group is complete}
\author{Vladimir Tolstykh}
\address{Vladimir Tolstykh\\ Department of Mathematics\\ Yeditepe University\\
34755 Kay\i\c sda\u g\i \\
Istanbul\\
Turkey}
\email{vtolstykh@yeditepe.edu.tr}
\maketitle

\begin{abstract}
Dyer and Formanek (1976) proved that if $N$ is a free nilpotent
group of class two and of rank $\ne 1,3,$ then the automorphism
group $\aut N$ of $N$ is complete. The main result of this paper states
that the automorphism group of an infinitely generated free
nilpotent group of class two
is also complete.
\end{abstract}

\section*{\it Introduction}

According to the result by J.~Dyer and E.~Formanek
\cite{DFo}, the automorphism group of a finitely
generated free two-step nilpotent group is complete
except in the case when this group is a one- or
three-generator (the three-generator groups have
automorphism tower of height two).  The purpose of
this paper is to prove that the automorphism group of
an infinitely generated free two-step nilpotent group
is also complete.  (Recall that a group $G$ is
said to be {\it complete} if $G$ is centreless
and every automorphism is inner.)

The paper may be considered as a contribution to
the study of automorphism towers of relatively free
groups. The study arose from conjectures of
G.~Baumslag that automorphism towers of finitely
generated absolutely free and free nilpotent groups
must be very short (strictly speaking, his conjecture
in the nilpotent case concerns finitely generated
torsion-free nilpotent groups \cite[problem
4.9]{Kou}).  The case of absolutely free groups was
considered in the paper \cite{DFo_abs} by Dyer and
Formanek:  they proved the automorphism group of a
finitely generated non-abelian free group $F$ is
complete, that is $\aut{\aut F} \cong \aut F$ (in
fact, the automorphism group of any non-abelian free group
is complete \cite{To_abs}). In \cite{DFo_app} Dyer
and Formanek obtained the following generalization
of the results from \cite{DFo_abs}: the automorphism
group of any group $F/R',$ where $F$ is absolutely
free of finite rank and $R$ is a characteristic
subgroup of $F$ lying in the commutator
subgroup $F'$ of $F,$
is complete provided that
$F/R$ is approximated by torsion-free nilpotent
groups. In particular, any finitely generated non-abelian
free solvable group has complete automorphism group.
The above cited paper
\cite{DFo} and the present paper give full description
of automorphism towers of free two-step nilpotent
groups.

Let $N$ denote an infinitely generated free two-step
nilpotent group. The ideas we use in the present paper
are closely related to those of Dyer and Formanek
\cite{DFo}. Thus, like the cited authors, we prove in
the last section that after multiplication by a
suitable inner automorphism of $\aut N$ any
automorphism of $\aut N$ preserves the elements of the
subgroup $\inn N$ and a fixed automorphism of $N$
which inverts all members of some basis of $N$ (we
call such automorphisms of $N$ {\it symmetries}; note
also that in \cite{DFo} a similar result, in
a part concerning the symmetry, is formulated a bit
weaker). However, instead of further analysis of the
action of the transformed automorphism, say, $\Delta$ of
$\aut N$ on {\it all generators} of $\aut N,$ we prefer to
prove that $\Delta$ preserves all
$\iat$-automorphisms, and hence the elements of the
conjugacy class of {\it all symmetries} (Theorem
\ref{DeltaFixesIAs}). This enables us to prove that
$\Delta$ preserves all elements of $\aut N.$
(Generally speaking, if $K$ is a conjugacy class of a
group $G$ such that the centralizer of $K$ in $G$ is
trivial, then any automorphism of $G$ which fixes
all elements of $K$ necessarily fixes all elements of $G.$)

Statements which are formulated similarly or exactly
the same as some statements from \cite{DFo} can be
also found in Sections \ref{invs} and \ref{conjs&symms}; all
such statements, mostly with different proofs, will be
specially indicated in the main body of the paper. The proof of
the completeness of $\aut N$ given in this paper
follows, nevertheless, an alternative general plan.

The specific feature of our proofs is the method we
usually use to show invariance of a subset of $\aut N$
under automorphisms of the group. The method is based
on the following general observation that came from
model theory: if a subset of an algebraic structure is
the set of realizations of a formula of a certain
logic with parameters in the structure, or, in
model-theoretic terms, if this subset can be {\it defined} by the
mentioned formula, then every automorphism of the
structure, which fixes each of the parameters, setwise
fixes the subset. In particular, if a subset is
definable by a formula without parameters then the
subset is invariant under all automorphisms of the
structure. In such situations algebraists used to say
that the subset of the structure can be characterized
in terms of basic operations. Thus, the reader who is
not familiar with model-theoretic terminology can
substitute his or her own arguments where necessary.
We should stress, however, that this paper does not
assume familiarity with model theory.

Usually, to define subsets we shall use formulae of
the first-order logic or the monadic second-order
logic (which allows quantification by arbitrary subsets
of a structure).  A subset of a structure definable by
means of first-order logic is simply called a {\it
definable} subset. The use of monadic second-order
is not actually particularly deep: we are just trying
to express the fact that characterization of some subsets
in $\aut N$ requires higher-order relations.

After preliminary Section \ref{basics} outlining
terminology and background material,
we begin by showing that a family of all
involutions of $N$ which are symmetries modulo
the subgroup $\ia N$ is definable in $\aut N$ (Lemma
\ref{SymmsBasics}). Then we prove that the subgroup
$\ia N$ itself is a definable, and hence a
characteristic subgroup of $\aut N$ (Proposition
\ref{IAsAreDef}).  In the same section we prove definability
in $\aut N$ modulo $\ia N$ for one more family of
involutions of $N,$ for extremal involutions (we call
an automorphism $\f$ of $N$ an {\it extremal}
involution, if there is a basis of $N$ such that $\f$
inverts some element of this basis element and fixes
others; the term is chosen in analogy with classical
group theory).

In Section \ref{conjs&symms} we prove definability of
conjugations and symmetries in $\aut N.$
First, using involutions
extremal modulo $\ia N$ we prove that the set of all
conjugations by powers of primitive elements is
definable in $\aut N$ (Lemma \ref{Conjs-by-PPE}; note
that a similar result holds for the automorphism
groups of non-abelian free groups \cite{To_abs}).
Lemma \ref{Conjs-by-PPE} implies that the subgroup
$\inn N$ is characteristic in $\aut N.$ Next, by means
of monadic second-order logic we define in $\aut N$
symmetries; this involves symmetries modulo $\ia N$
and normalizers of free generating sets of the (free
abelian) group $\inn N$ (Lemma
\ref{Only-Conjs-of-a-Symm}).

The main result of Section \ref{ia-stabs}, the next to
the last in this paper, states that the subgroup
$\iat_{\tau}(N)$ of $\aut N$ consisting of all
$\iat$-automorphisms which stabilize a given primitive
element $x$ of $N$ (in fact, any element in $x N',$
where $N'$ is the commutator subgroup of $N$) is
definable in $\aut N$ by means of monadic second-order
logic with the parameter $\tau,$ where $\tau$ is
conjugation by $x$ (Theorem \ref{Def-of-IAtau}). The
stabilizers $\iat_{\tau}(N)$ are involved in a second-order
modelling in $\aut N$ the primitive elements of $N$
and play a crucial role in the proofs in the last
section, which was briefly described above.

The author would like to express his gratitude to
Professor Oleg Belegradek for friendly attention
to this research and to thank the referee for
helpful comments and suggestions.

\section{\it Basic concepts and notation} \label{basics}

Everywhere in this paper $N$ denotes an infinitely
generated free two-step nilpotent group, $N'$ stands
for the commutator subgroup of $N$ and $A$ for the
free abelian group $N/N'.$ We denote by $\bar{\phantom a}$
the natural homomorphism $N \to A,$ and use
the same symbol to denote the corresponding induced
homomorphism $\aut N \to \aut A.$

For any two-step nilpotent group the commutator
subgroup is contained in the centre of this group; for
a free two-step nilpotent group the centre, a
free abelian group, is exactly
the commutator subgroup (\cite[5.7]{MKS}, \cite[ch. 3, \S \ 1]{HN}). We denote by $\tau_a$ the inner
automorphism of $N,$ or conjugation, determined
by an element $a \in N.$ Since $N'$ is the
centre of the group $N,$ then $\tau_a=\tau_b$
if and only if $a \equiv b (\mod N').$ Hence $\inn N,$
the group of all inner automorphisms of $N,$ is
isomorphic to $N/N',$ and, in particular, is a free
abelian group.

\begin{Th} \label{Maltsev's}
A set $\{x_i : i \in I\}$ is
a basis {\em(}free generating set{\em)} of $N$ if
and only if the set $\{\av x_i : i \in I\}$
is a basis of the free abelian group $N/N'.$
\end{Th}

Recall that the sufficiency part of the Theorem can be
proved by using the following two results: (1) if a
set $X \cup G'$ generates a nilpotent group $G,$ then
$G$ is generated only by $X$ itself \cite[Cor.
10$\cdot$3$\cdot$3]{MHall} and (2) if $X$ generates a free
nilpotent group $G$ and $\av X$ is a basis of $G/G',$
then $X$ is a basis of $G$ \cite[\S \ 4]{Malt}
(see also \cite[ch 3, \S \ 1, ch. 4, \S \ 2]{HN}).

\begin{Cor} \label{AvIsSurj}
Every automorphism of $A$ is induced by an
automorphism of~$N.$
\end{Cor}

The kernel of the induced homomorphism $\aut N \to
\aut A,$ the subgroup of $\iat$-{\it automorphisms},
is denoted as usual by $\ia N.$ Considering the action
of the elements of $\ia N$ on a fixed basis of $N,$ one
sees that $\ia N$ is isomorphic to an infinite
Cartesian power of $N',$ and therefore is a
torsion-free abelian, but not free abelian group as in
the case of finitely generated free two-step nilpotent
groups \cite{Bach,DFo}.

We shall work with involutions in the automorphism
group of the abelianization $A.$ Suppose that $f$ is
an involution of the group $A.$ Write $A^+_f$ (in the
manner of R.~Baer) for the fixed-point subgroup of $f$
and $A^-_f$ for the subgroup of elements $\{a\}$ such
that $f a=-a.$ An involution $f \in \aut A$ is
diagonalizable in some basis of $A$ if and only if
$$
A=A^+_f\oplus A^-_f
$$
(note that the latter property does not hold for {\it
all} involutions in $\aut A;$ see Theorem \ref{HRCanForms}
below). We call a diagonalizable involution $f$ a
$\vk$-{\it involution}, where $\vk$ is a cardinal, if
$$
\vk=\min(\rank A^+_f,\rank A^-_f).
$$
A standard argument proves that

\begin{Lem} \label{SoftComm}
Diagonalizable involutions $f$ and $g$ from $\aut A$
commute if and only if
$$
A=(A^+_f \cap A^+_g) \oplus  (A^+_f \cap A^-_g) \oplus
(A^-_f \cap A^+_g) \oplus
(A^-_f \cap A^-_g).
$$
\end{Lem}

In general, the structure of involutions in $\aut A$
is described by the following

\begin{Th} \label{HRCanForms}
Any involution $f$ in the group $\aut A$ has
a basis $B$ of $A$ such that $f b=\pm b$ or
$f b \in B$ for each $b \in B.$
\end{Th}

The result is essentially known for free abelian
groups of finite rank (it immediately follows from
Lemma 1 in \cite{HuaRei} by L.~K.~Hua and I.~Reiner);
we give a sketch of the proof in the case of infinite
rank.

{\it Proof of Theorem} \ref{HRCanForms}. The
fixed-point subgroup $A^+_f$ is a free summand of $A.$
Let $R$ be a direct complement $\Fix(\f)$ to $A$ and
let $\{r_i : i \in I\}$ be a basis of $R.$
For every $i \in I$ the element $f r_i +r_i$ is in
$\Fix(\f).$ Suppose that
$$
J=\{i \in I : f r_i +r_i \in 2A\} \text{ and } K=I\setminus J.
$$
For every $j \in J$ we have $f r_j+r_j=2t_j,$
where $t_j \in \Fix(\f);$ therefore if $s_j=r_j-t_j,$
then $\f(s_j)=-s_j.$

The index set $K$ can be considered as a well-ordered.
One can then construct by transfinite induction
a family $\{s_k : k \in K\}$ of elements of $A$ such that
\begin{itemize}
\item[(i)] the set $\{s_k : k \in K\} \cup \{s_j : j \in
J\}$ generates a direct complement of $\Fix(A)$
to $A;$

\item[(ii)] $f s_k=-s_k+u_k,$ where $u_k=0$ or
$u_k$ is a unimodular element of $A$ (that is,
can be found in some basis of $A$) and
\item[(iii)] the family consisting of
all non-zero elements $u_k$ forms
a basis of a direct summand of $\Fix(\f).$
\end{itemize}

Let $I_1=\{i \in I : f s_i=-s_i\}$ and $I_2 =I\setminus I_1.$
By the construction the subgroup $\str{f s_i+s_i : i \in I_2}$
is a direct summand of $\Fix(\f)$:
$$
\Fix(\f)= \str{f s_i+s_i : i \in I_2} \oplus V.
$$
Suppose that $\cV$ is a basis of $V.$ Then the following basis of $A$
$$
\{s_i : i \in I_1\} \cup \{s_i,f s_i : i \in I_2\} \cup \cV
$$
satisfies the conditions of the Theorem.

\section{\it Coming down to the abelianization} \label{invs}

We shall use throughout the paper two conjugacy
classes of involutions in the group $\aut N.$ The first
conjugacy class consists of involutions similar to
involutions one can find in standard generating sets
of automorphism groups of finitely generated two-step
free nilpotent group. We shall call $\f
\in \aut N$ an {\it extremal} involution if there is a basis of $N$
such that $\f$ inverts exactly one element of this basis and
fixes other elements; we shall also call any basis of
$N$ on which $\f$ acts in such a way a {\it canonical}
basis for $\f.$ Any {\it symmetry} $\theta,$ a member of the second
conjugacy class, has a basis of $N$ such that $\theta$ takes each
element of this basis to the inverse (a {\it canonical}
basis for $\theta$).

Note that both described classes are used in the cited
paper \cite{DFo} by Dyer and Formanek. In this section
we obtain a first-order characterization of both
symmetries and extremal involutions in $\aut N$ modulo
$\ia N.$

\begin{Lem} \label{SymmsBasics}
\mbox{\em (a)} Let $\theta$ be a symmetry. Then
for every $\iat$-automorphism $\alpha$
$$
\theta \alpha \theta=\alpha\inv;
$$
in particular, the automorphism $\theta\alpha$ is
an involution in $\aut N.$


{\em (b)} Let $\theta$ be a symmetry, $\sigma \in \aut
N$ and $\alpha \in \ia N.$ If both $\sigma$ and
$\s\alpha$ commute with $\theta,$ then
$\alpha=\operatorname{id}.$

{\em (c)} An involution $\theta \in \aut N$
is a symmetry modulo $\ia N$ {\em(}that is, has
the form $\theta^* \beta$ for some symmetry
$\theta^*$ and $\beta \in \ia N${\em)} if and only if
\begin{equation} \label{ThreeConjs}
\text{a product
of any three conjugates of $\theta$ is an involution.}
\end{equation}

{\em (d)} The family of all involutions which
are symmetries modulo $\ia N$ is definable
in $\aut N.$
\end{Lem}

\begin{proof}
(a) Since $\theta$ inverts all elements of some basis
$\cB$ of $N,$ then $\theta$ inverts modulo $N'$ all
elements of $N.$ This implies that $\theta$ preserves
any commutator $[a,b] \in N'$:
$$
\theta [a,b] =[\theta a, \theta b]=[a\inv,b\inv]=[a,b\inv]\inv=
([a,b]\inv)\inv=[a,b],
$$
and hence all elements of $N'.$

Let $x$ be an element of $\cB$ and $\a x=x c,$
where $c \in N'.$ Then we have that
$$
\theta \alpha \theta x=\theta \alpha (x\inv)=
\theta (x\inv c\inv)=x c\inv=\alpha\inv x.
$$

(b) The conditions $\theta \sigma \theta=\sigma$
and $\theta \sigma \alpha \theta=\sigma\alpha$
imply by (a) that $\alpha\inv=\alpha,$ whence
$\alpha=\id,$ because $\ia N$ is a torsion-free group.

(c) According to (a), any symmetry modulo $\ia N$ is an
involution. Further, an involution from
$\aut N$ is a symmetry modulo $\ia N$ if and only if
its image under the induced homomorphism
$\aut N \to \aut A$ (recall that $A=N/N'$) is equal to
$-\id_A.$ Then a product of any three (conjugate) symmetries
modulo $\ia N$ is again a symmetry modulo $\ia N,$
and therefore an involution.

Let us prove the converse. It follows
from Corollary \ref{AvIsSurj} that if $\s \in \aut N,$
then any conjugate of $\av \s$ in $\aut A$ can be
lifted to a conjugate of $\s$ in $\aut N$:  if $s'
\sim \av \s$ in $\aut A,$ then there is $\s' \in \aut
N$ such that $\s' \sim \s$  and $\av \s' =s'$ ($\sim$
denotes the conjugacy relation).
Then it suffices to
prove that $-\id,$ or a unique non-trivial central
element of $\aut A,$ is the only involution with
the property \eqref{ThreeConjs} in the group $\aut A.$

We shall base our argument on Theorem
\ref{HRCanForms}. It is quite clear in view of this
Theorem that the latter statement follows from the
fact that the similar result holds for the
automorphism group of a two-generator free abelian
group, or, equivalently, for the group $\gl{2,\Z}.$
Indeed, if $f$ is an involution in $\aut A$ and
$B$ is a basis of $A$ with the property
described in Theorem \ref{HRCanForms}, then
there are two elements $u,v$ of $B$ such
that both subgroups $\str{u,v}$ and
$\str{B \setminus \{u,v\}}$ are $f$-invariant.
Thus, one can make conjugates of $f$ changing
its action on $\str{u,v},$ but preserving
the action on $\str{B \setminus \{u,v\}}.$

Theorem \ref{HRCanForms} implies also that any matrix of
order two from $\gl{2,\Z}$ which is not in the centre of
this group is conjugate to the matrix
$$
X=
\begin{pmatrix}
1 & 0 \\
0 & -1
\end{pmatrix}
$$
or to the matrix
$$
Y=
\begin{pmatrix}
0 & 1 \\
1 & 0
\end{pmatrix}.
$$
One easily checks that neither the product of the
following three conjugates of $X$
$$
X_1=
\begin{pmatrix}
1 & 0              \\
0 & -1
\end{pmatrix},\quad
X_2=
\begin{pmatrix}
1 & 0         \\
2 & -1
\end{pmatrix},\quad
X_3=
\begin{pmatrix}
-1 & 2 \\
0 & 1
\end{pmatrix},
$$
nor the product of the following conjugates of $Y$
$$
Y_1=
\begin{pmatrix}
1 & 0           \\
1 & -1
\end{pmatrix},\quad
Y_2=
\begin{pmatrix}
1 & 0             \\
-1 & -1
\end{pmatrix},\quad
Y_3=
\begin{pmatrix}
0 & 1 \\
1 & 0
\end{pmatrix}
$$
has order two.

(d) By (c).
\end{proof}

As an immediate consequence of the previous Lemma
we have:

\begin{Prop} \label{IAsAreDef}
\mbox{\em (a)} An automorphism $\s \in \aut N$ is
in the subgroup $\ia N$ if and only
if $\s$ can be written as a product
of two involutions with \eqref{ThreeConjs}.

{\em (b)} $\ia N$ is a definable subgroup
of $\aut N.$
\end{Prop}

Let us obtain one more consequence.

\begin{Prop} \label{Centreless}
The group $\aut N$ is centreless.
\end{Prop}

\begin{proof}
Since the induced homomorphism $\aut N \to \aut A$
is surjective, then an element $\s$ from the centre
of $\aut N$ should induce a central element of
the group $\aut A.$ The only elements in the
centre of $\aut A$ are $\pm \id_A,$ and therefore
$\av\s=\pm \id_A.$ In the case when $\av\s=\id_A$
we have $\s \in \ia A,$ but there is no non-trivial
$\iat$-automorphism which commutes with symmetries
(Lemma \ref{SymmsBasics} (a)). Assuming $\av\s=-\id_A$
we obtain that $\s$ is a symmetry modulo $\ia N,$ that
is, has the form $\theta \beta$ for some symmetry
$\theta$ and an $\iat$-automorphism $\beta.$ But
again by Lemma \ref{SymmsBasics} we have for every
$\iat$-automorphism $\alpha$ that
$$
\theta \beta \alpha \beta\inv \theta=(\text{$\ia N$ is an abelian group})=
\theta \alpha \theta=\alpha\inv,
$$
and the result follows.
\end{proof}

Next is a first-order characterization of extremal
involutions modulo $\ia N$ in the group $\aut N.$ The
fact that $\ia N$ is a definable subgroup of $\aut N$
(Proposition \ref{IAsAreDef}) enables us to involve
the structures related to $\ia N$ and the subgroup
$\ia N$ itself into first-order characterisations
of subsets of $\aut N.$

\begin{Prop}
Let $K(f)$ denote the conjugacy class of an
automorphism $f$ of $A$ in the group $\aut A.$
An involution $\f \in \aut N$ is an extremal
modulo $\ia N$ if and only if
\begin{itemize}
\item[(i)] $\av \f$ is not a square in $\aut A;$
\item[(ii)] the set $K^2(\av\f)=K(\av\f)K(\av\f)$ contains
no elements of order three and
\item[(iii)] all involutions in $K^2(\av\f)$
are conjugate.
\end{itemize}
The family of involutions extremal modulo $\ia N$
is definable in $\aut N.$
\end{Prop}

\begin{proof} If $\f$ is an extremal involution modulo $\ia N,$ then
there is a basis $\{x\} \cup Y$ of $N$ such that
\begin{alignat*} 2
\f x &\equiv x\inv &&(\mod N'),\\
\f y &\equiv y     &&(\mod N'),\quad \forall y \in Y.
\end{alignat*}
Let $f=\av\f.$ Then $\pm \av x$ are the
only unimodular elements in $A^-_f=\{a : f a=-a\}.$
Suppose that $f=g^2,$ where $g \in \aut A.$
We then have $f(g\av x)=-g\av x,$ since
$g$ commutes with $f.$
Hence $g\av x=\pm \av x$ and the equation $f=g^2$
is impossible.

Let us check (ii). Consider a natural homomorphism
$\widehat{\phantom{\sigma}}$ from $A$ onto $A/2A,$ the quotient
group of $A$ by the subgroup of even elements;
let $\widehat{\phantom{\sigma}}$ denote also the corresponding
induced homomorphism $\aut A \to \aut{A/2A}.$ Take
an automorphism $s \in \aut A$ of order
three. We claim that the image of
$s$ under the homomorphism $\widehat{\phantom{\sigma}}$ is non-trivial. This
will imply that a product of any two conjugates
of $\av\f,$ an element of $K^2(\av\f),$
cannot have order three, since
the image of $\av\f$ in $\aut{A/2A}$
under $\widehat{\phantom{\sigma}}$ is trivial.

The kernel of the endomorphism $s^2+s+\id$ of the
group $A$ has non-trivial elements modulo $2A$ (it
contains, in particular, the subgroup $(s-\id)A).$ Let
$a$ be such an element. Suppose that $\widehat s=\id.$
Hence we have
$$
0=\widehat s^2 \widehat a+\widehat s\,\widehat a+\widehat a=
3\widehat a=\widehat a,
$$
a contradiction.

To check the last condition (iii) we use Lemma \ref{SoftComm}.
The image of $\f$ in $\aut A$ is a 1-involution (see
Section \ref{basics} for definitions), and
according to Lemma \ref{SoftComm} a product of
any two conjugate and commuting 1-involutions is a 2-involution whose
fixed point-subgroup has rank equal to $\rank A$
(recall that $\rank A$ is an infinite cardinal).
Thus, any two involutions from  $K^2(\av\f)$ are
conjugate.

Conversely, suppose that $\av\f,$ where $\f$ is
an involution from $\aut N,$ satisfies the conditions
(i-iii). The condition (ii) imply that $f=\av\f$ is
diagonalizable. Indeed, suppose, towards a contradiction,
that there is a basis $B$ of $A$ such that
$f B \subseteq \pm B$ and there exist two elements
$u,v$ of $B$ taken by $f$ to one another (Theorem \ref{HRCanForms}). It is
easy then to find conjugates $f',f''$ of $f$ such that
$f'f''$ has order three: we require that
$f' b=f'' b=f b$ for each $b \in B\setminus \{u,v\}$
and
$$
\begin{cases}
f' u=-v,\\
f'v=-u,
\end{cases}
\quad
\begin{cases}
f'' u=u+v,\\
f'' v=-v.
\end{cases}
$$

Thus, we have that $f$ is a diagonalizable, and hence
a $\vk$-involution for some cardinal $\vk.$ If $\vk >
1,$ then one finds in $K^2(f)$ not only 2-involutions,
but also, for instance, 4-involutions and (iii) fails.

To complete the proof we have to show that a
1-involution $g$ such that $\rank A^-_g=\rank A$ is a
square in $\aut A.$ Assume $\{b\} \cup C$ is a basis
of $A$ such that $g b=b$ and $g c=-c$ for each $c\in
C.$ Then we construct $h \in \aut A$ whose square
is $g$ putting $h b=b$ and having in mind a well-known
relation
$$
\begin{pmatrix}
0 & -1 \\
1 & 0
\end{pmatrix}^2
=
\begin{pmatrix}
-1 & 0\\
0 & -1
\end{pmatrix}.
$$
in order to define the action of $h$ on $C.$
This completes the proof of the Proposition.
\end{proof}

We close this section with a simple result on the
conjugation action of extremal involutions on $\ia N.$

\begin{Lem} \label{PMs-of-an-Ext}
Let $\f$ be an extremal involution and $\cB=\{x\} \cup
Y$ a canonical basis for $\f$ {\em(}that is $\f x=x\inv$
and $\f y=y$ for each $y \in Y${\em).} Suppose that
$\alpha$ is an $\iat$-automorphism. Then

{\em (a)} $\f \alpha \f=\alpha$ if
and only if $\alpha$ moves the elements of $\cB$
as follows:
\begin{alignat*}2
&\alpha x &&= x[x,t],\\
&\alpha y &&=y d_y, \quad y \in Y, \nonumber
\end{alignat*}
where $t$ is an element of the
subgroup $\str Y$ generated by $Y$
and $d_y \in \str Y{}'$ for each
$y \in Y;$

{\em (b)} $\f \alpha \f=\alpha\inv$ if
and only if $\alpha$ moves the elements of $\cB$
as follows:
\begin{alignat*}2
&\alpha x &&= x c,\\
&\alpha y &&=y [x,a_y], \quad y \in Y, \nonumber
\end{alignat*}
where $c \in \str Y{}'$ and $a_y \in \str Y$
for each $y \in Y.$
\end{Lem}

\begin{proof}
The important point is that the restriction of $\f$ on
$N'$ is a diagonalizable involution, and $[x,N]$ and
$\str Y{}'$ are its $(-)$- and $(+)$-subgroup,
respectively. Really, $\f$ fixes or inverts natural
generators of $N',$ basis commutators $[u,v],$ where
$u,v$ are distinct elements of $\cB.$ Clearly $\f
[u,v]=[u,v]$ if neither $u=x,$ nor $v=x.$ On the other
hand, if $v \ne x,$ then $\f [x,v]=[x,v]\inv.$

Let us prove now, for example, (a). Assume that $\alpha x=x d,$
where $d \in N'.$ We have
$$
\f \alpha \f x=\alpha x \iff
x \f(d\inv)=xd \iff \f d=d\inv \iff d \in [x,N].
$$
Similarly, if $\alpha y=y d_y,$ then
$$
\f \alpha \f y=\alpha y \iff \f d_y =d_y \iff d_y \in \str Y{}',
$$
as required.
\end{proof}

{\sc Remark.} Our starting point, the proposition stating
that $\ia N$ is a characteristic subgroup of $\aut
N,$ is the same as in the paper \cite{DFo} by Dyer
and Formanek. They prove (the proof of Lemma 3 on
page 273) that $\ia{N_k}$ is the Hirsch-Plotkin
radical of $\aut{N_k},$ that is the maximal locally
nilpotent subgroup of $\aut{N_k},$ where $N_k$ is a
$k$-generator ($k < \infty$) free two-step nilpotent
group; therefore $\ia{N_k}$ is characteristic in
$\aut{N_k}.$  On the same page of \cite{DFo} one finds
a statement on the conjugation action of a symmetry on
$\ia{N_k}$ and the proof that $\aut{N_k}$ is
centreless. These facts correspond to Lemma
\ref{SymmsBasics} (a) and Proposition
\ref{Centreless}, respectively.

\section{\it Characterizing conjugations and symmetries} \label{conjs&symms}

We begin with the proof of definability of
conjugations by powers of primitive elements of $N$ in
$\aut N$ (recall that an element of the group $N$ is
said to be {\it primitive} if it is a member of some
basis of this group). The mentioned conjugations
generate the subgroup $\inn N,$ and hence the latter
is characteristic in $\aut N.$ The rest of the section
is devoted to a characterization of symmetries and
consideration of the ways in which the basis sets of
$N$ and the primitive elements of $N$ can be modelled
in $\aut N.$

\begin{Lem} \label{Conjs-by-PPE}
An automorphism $\alpha \in \ia N$ is conjugation
by a power of a primitive element of $N$ if
and only if there is an involution $\f$
extremal modulo $\ia N$ such that
\begin{itemize}
\item[(i)] $\f \alpha \f=\alpha\inv;$
\item[(ii)] if $\f$ and $\psi,$ where $\psi$
is an extremal involution modulo $\ia N,$ are distinct
and commuting modulo $\ia N,$ then $\psi$
commutes with $\alpha.$
\end{itemize}
The family of conjugations
by powers of primitive elements is definable
in $\aut N.$
\end{Lem}

\begin{proof}
If $\tau$ is conjugation by, say, an $m$th power
a given primitive element $x \in N,$ then one
easily finds an extremal involution
$\f$ which inverts $x,$ whence $\f \tau \f =\tau\inv.$
If further $\psi$ is an extremal involution modulo $\ia N$
satisfying the conditions of the Lemma,
then by Lemma \ref{SoftComm} there exist
a basis of the free abelian group $A$ in which
both $\avst\f$ and $\av\psi$ are diagonalizable.
Since $\avst\f \ne \av\psi,$ then $\av\psi\, \avst x=\avst x.$
This implies that $\psi x =x c,$ where $c \in N',$
and hence
$$
\psi \tau \psi =\psi \tau_x^m \psi=\tau_{x c}^m =\tau_x^m =\tau.
$$

Let us prove the converse. The conditions (i) and (ii) deal with
the conjugation action on $\ia N$ and the commutativity
modulo $\ia N;$ this then allow us to
assume that $\f$ is an extremal involution.
Let $\cB=\{x\} \cup Y$ be a canonical
basis of $N$ for $\f.$ According to Lemma \ref{PMs-of-an-Ext},
an IA-automorphism $\alpha$ with (i) moves the elements of $\cB$ as follows:
\begin{alignat*}2
&\alpha x &&= x c_x,\\
&\alpha y &&=y [x,a_y], \quad y \in Y, \nonumber
\end{alignat*}
where $c_x \in \str Y{}'$ and $a_y \in \str Y$
for each $y \in Y.$

Any extremal involution $\f_y$ which
inverts $y \in Y$ and
whose action on $\cB$ is canonical commutes
with $\f.$ Therefore if (ii) holds,
then, again by Lemma \ref{PMs-of-an-Ext}, $\f_y c_x=c_x$
and $\f_y [x,a_y]=[x,a_y]\inv$ for each $y \in Y.$
As for $c_x,$ we have that
$$
c_x \in \bigcap_{y \in Y} \str{\cB \setminus \{y\}}',
$$
and this, along with $c_x \in \str{\cB \setminus
\{x\}}',$ implies that $c_x=1.$ For every $y \in Y$ the
commutator $[x,a_y]$ must be equal to $[y,b_y]$ for
some $b_y$ in $\str{\cB \setminus \{y\}}$:
$[x,a_y]=[y,b_y].$ It follows that
$[x,a_y]=[x,y]^{k_y}$ for a suitable integer $k_y \in
\Z.$ So the action of $\alpha$ on $\cB$ looks like
\begin{align*}
\alpha x &=x,\\
\alpha y &=x^{k_y} y x^{-k_y}.
\end{align*}
To prove that $\alpha$ is conjugation we have to prove
that the integers $k_y$ are the same.

Choose an element $z \in Y$ and consider a basis
$$
\cB'=\{x\} \cup \{z\} \cup \{yz : y \in Y, y \ne z\}
$$
of the group $N.$ The action of $\f$ on $\cB'$ is also canonical
and by applying the above arguments one can
deduce that
$$
\alpha(yz)=x^{m_y} (yz) x^{-m_y}
$$
for each $y \in Y \setminus \{z\}.$
Thus, for every $y \in Y\setminus \{z\}$
we have $m_y=k_y=k_z,$ and $\alpha$ is conjugation
by a power of $x,$ as desired.
\end{proof}

\begin{Cor}
The subgroup of all conjugations $\inn N$ is definable
in $\aut N.$
\end{Cor}

\begin{proof}
Every element of an infinitely generated
free abelian group (in particular, every
element of $\inn N$) can be written
as a product of two unimodular (primitive)
elements.
\end{proof}

{\sc Remark.} It can been seen quite easily that every
element of a free abelian group is actually a product
of at most three unimodular elements \cite{Smir}.

The fact that a set $B$ is a basis of a free abelian
group $\str{G,+}$ can be easily expressed by a formula
of monadic second-order logic. This formula
may be chosen as a `translation' of the
following statement: a subset $B$ of
$G$ is a basis of $G$ if and only if
$$
G=\str b \oplus \str{B \setminus \{b\}}
$$
for each $b \in B.$

We shall call a basis $B$ of the free abelian group
$\inn N$ a {\it basis set} of conjugations. By Theorem
\ref{Maltsev's} there is a basis $\cB$ of $N$ such
that conjugations in $B$ are determined by the
elements of $\cB$: $B=\{\tau_b : b \in \cB\}.$
Consider a symmetry $\theta^*$ which inverts all
elements of $\cB$ and let $\operatorname{N}(B)$ denote
the normalizer of $B$ in $\aut N.$

\begin{Lem} \label{Only-Conjs-of-a-Symm}
Let $\theta$ be a symmetry modulo $\ia N.$ Then the
following statements are equivalent:
\begin{itemize}
\item[(i)] $\theta$ has the form $\theta^* \alpha^2,$
where $\alpha \in \ia N;$

\item[(ii)] $\theta$ commutes modulo the subgroup
$\iat^2(N)$ with each element of $\operatorname{N}(B),$
where $\iat^2(N)$ is the subgroup of $\ia N$
generated by squares of the elements
of $\ia N.$

\end{itemize}

\end{Lem}

One immediately deduces from the Lemma that

\begin{Prop} \label{DefOfSymms}
There is a monadic second-order formula
which is satisfied in $\aut N$ exactly by
symmetries.
\end{Prop}

\begin{proof}
Any involution of the form $\theta^* \alpha^2$ is a
symmetry, namely a conjugate of the symmetry
$\theta^*$: $\theta^* \alpha^2 = \alpha\inv \theta^*
\alpha$ (Lemma \ref{SymmsBasics} (a)).
\end{proof}

{\sc Remark.} We use an idea Dyer and Formanek use in
\cite{DFo} in order to identify the image of a chosen
symmetry under an automorphism of $\aut{N_k},$ where
$N_k$ is a free two-step nilpotent of finite rank $k.$
In the case of infinite rank this idea can be applied,
however, with an optimal effect, since one finds among
the realizations in $\aut{N_k}$
of the condition (ii) from Lemma \ref{Only-Conjs-of-a-Symm}
involutions which are not necessarily symmetries (for
example, if $\{x_1,\ldots,x_k\}$ is a basis of $N_k$
corresponding to our basis $\cB,$ then an involution
$\theta^* \tau_{x_1} \ldots \tau_{x_n}$ is such a
realization).

{\it Proof of Lemma} \ref{Only-Conjs-of-a-Symm}.
Every element in the normalizer
of $B$ in $\aut N$ can be written in
the form $\pi \beta,$ where $\pi$ acts
on the basis $\cB$ as a permutation
(and hence commutes with $\theta^*$) and
$\beta \in \ia N.$ Therefore
\begin{align*}
\pi \beta (\theta^* \alpha^2) \beta\inv \pi\inv &=
\pi \theta^* (\alpha \beta\inv)^2 \pi\inv \\
&=  \theta^* \pi (\alpha \beta\inv)^2 \pi\inv \equiv
\theta^* \alpha^2 (\mod \iat^2(N)).
\end{align*}

Conversely, preserving notation we have just
introduced, suppose that for any $\pi$ and $\beta$
$$
\pi \beta (\theta^* \gamma) \beta\inv \pi\inv \equiv \theta^* \gamma (\mod \iat^2 (N)),
$$
where $\gamma \in \ia N.$ It then follows that
for any $\pi$
\begin{equation}
\pi \gamma \pi\inv \equiv \gamma (\mod \iat^2(N)).
\end{equation}
We claim that $\gamma$ is a square in $\ia N.$

The assumption of the existence of $b \in \cB$
such that $\gamma b =b d_b^2,$ where $d_b \in N'$
trivially guarantees the conclusion. Really,
let $t \in \cB \setminus \{b\}$ and
$\gamma t=t c_t,$ where $c_t \in N'.$ Taking
$\pi$ such that $\pi b=t$ we have
$$
\gamma t \equiv \pi \gamma \pi\inv t (\mod N'{}^2) \Rightarrow
t c_t \equiv t \pi(d_b^2) (\mod N'{}^2) \Rightarrow c_t \equiv 1 (\mod N'{}^2).
$$

Therefore, we assume that for each $b \in \cB$
$\gamma b=b c_b,$ where $c_b$ is an element
of $N'$ which is not a square in $N'.$
Take an arbitrary $b \in \cB$ and suppose
that the word $c_b$ (in the letters $\cB$) has non-trivial
occurrences of an element $a \in \cB$:
$$
c_b =[b,a]^k [b,u_b] [a,v_b] d_b,
$$
where the words $u_b,v_b,$ and $d_b \in N'$
contain no occurrences of both $a$ and $b,$
and $v_b \notin N'$ or $k \ne 0.$
Write in analogous way the action of $\gamma$
on the element $a$:
$$
\gamma a = a [a,b]^m [a,u_a] [b,v_a] d_a.
$$
Assume that the elements $u_a,v_a,d_a,u_b,v_b,d_b$
are the words in the letters $b_1,\ldots,b_n \in \cB.$
Since $\cB$ is infinite, there is an automorphism $\pi \in \aut N$ which
preserves $\cB,$ takes $a$ and $b$ to each other,
and such that
$$
\pi \{b_1,\ldots,b_n \} \cap \{b_1,\ldots,b_n\} =\varnothing.
$$
By (\theequation) $\gamma b \equiv \pi \gamma \pi\inv b (\mod N'{}^2),$
and then
$$
[b,a]^k [b,u_b] [a,v_b] d_b \equiv [b,a]^m [b,\pi u_a] [a,\pi v_a] \pi d_a (\mod N'{}^2).
$$
This implies that the elements $u_a,v_a,u_b,v_b$
are all squares modulo $N'$ and $d_a,d_b \in N'{}^2.$
Hence $c_b \equiv [b,a]^k (\mod N'{}^2).$

Consider $t \in \cB \setminus \{b,a\}.$ We have
$\gamma t \equiv t[t,s]^n (\mod N'{}^2)$ for some $s \in \cB.$ To prove
that the element $[b,a]^k$ is a square, one may use
$\pi$ such that $\pi$ takes $t$ to $b$ and fixes all
elements in $\cB \setminus \{b,t\}$ (if $s \ne a$), or
$\pi$ which acts on $\{a,b,t\}$ as a cycle $(b,a,t)$
(if $s=a$).

The proof of Lemma \ref{Only-Conjs-of-a-Symm}
is now completed.
$\qed$

In view of Lemma \ref{Only-Conjs-of-a-Symm},
a symmetry $\theta$ satisfying the condition
(ii) from this Lemma
for a given basis set of conjugations $B$ will be called
{\it attached} to $B.$ One can then associate with a
pair $(B,\theta)$ a basis $\cC = \cB(B,\theta)$ of $N$
uniquely determined by the
following conditions:

\begin{itemize}
\item[(i)] $B=\{\tau_z : z \in \cC\};$
\item[(ii)] $\theta$ inverts all elements of $\cC$
(that is, acts on $\cC$ canonically in
terms introduced in Section \ref{invs}).
\end{itemize}
Indeed, if $\cB$ is a basis of $N$ which satisfies (i) and
$\theta_\cB$ is the symmetry whose action on $\cB$ is canonical, then
by Lemma \ref{Only-Conjs-of-a-Symm} $\theta=\alpha\inv
\theta_\cB \alpha$ for some $\iat$-automorphism
$\alpha.$ Hence $\theta$ inverts all elements
of the basis $\alpha\inv \cB$ of $N,$ and, moreover, this
basis satisfies (i). Finally, if $x$ is a
primitive element of $N$ (in particular,
an element of $\cB$) and $\theta (xc)=x\inv c\inv$
for some $c \in N',$ then $xc$ is the only
element of $x\,N'$ taken by $\theta$ to
the inverse: assuming $\theta (xd)=x\inv d\inv,$
where $d \in N',$ we have that
$$
\theta (xd)=\theta(x) \theta(d)=x\inv c^{-2} d;
$$
therefore $c^{-2} =d^{-2},$ or $c=d.$

In a more general setting {\it any} triplet
$(\tau,B,\theta),$ where $B$ is a basis set of
conjugations, $\theta$ a symmetry attached to $B$ and
$\tau \in B,$ codes a primitive element of $N$: we
assign to each such a triplet a unique element $x$ of the basis
$\cB(B,\theta)$ such that $\tau=\tau_x,$ or,
equivalently a unique element in $y N',$ where $y$ is
any primitive such that $\tau=\tau_y,$ taken by
$\theta$ to the inverse.

Two triplets $(\tau,B,\theta)$
and $(\tau',B',\theta')$ will code
the same primitive element of $N$ if and only
if the following conditions hold
\begin{itemize}
\item[(PE1)] $\tau=\tau'$ (and hence
there is a primitive $x \in N$ such
that $\tau=\tau'=\tau_x$);
\item[(PE2)] $\theta$ and $\theta'$ invert
the same element in $x N'.$
\end{itemize}

It is easy to see that if an $\iat$-automorphism
$\alpha$ preserves a primitive $y \in N,$ then
it preserves all elements in $y N'.$ Then
the condition (PE2) is equivalent
to the condition which states that
the $\iat$-automorphism $\theta\theta'$
fixes all elements in $x N'$ (this
will imply that $\theta x =\theta' x,$ and
hence (PE2) will hold).

Let $\tau$ be conjugation by a primitive
element, and let $\iat_\tau(N)$ denote
the subgroup of all $\iat$-automorphisms
which fix any $z \in N$ such that $\tau=\tau_z.$
Thus, in order to prove that the condition
(PE2) can be expressed in $\aut N$ by means of group
theory it suffices to obtain a characterization
of the subgroups $\iat_\tau(N).$ This is
the main subject of the next section.

\section{\it \mbox{\rm IA}-stabilizers} \label{ia-stabs}

\begin{Th} \label{Def-of-IAtau}
Let $B$ be a basis set of conjugations and $\theta$ a
symmetry attached to $B.$ Then for each $\tau \in B$
the subgroup $\iat_\tau(N)$ is definable with
the parameters $\tau,\theta,B$ in $\aut N$ by means of
monadic second-order logic.
\end{Th}

{\sc Remark.} In fact $\iat_\tau,$ where $\tau$ is
conjugation by a primitive element, is definable in
$\aut N$ by means of monadic second-order logic only
with the parameter $\tau.$ The use of other parameters
in the Theorem is a little more convenient for the
proofs in the next section.

\begin{proof}
As we saw in the previous section there is a unique
basis $\cB=\cB(B,\theta)$ of $N$ such that the set of
conjugations by elements of $\cB$ is $B$
and $\theta$ inverts all elements of $\cB.$
Assume $x$ is an element of $\cB$ such
that $\tau=\tau_x$ and write $\cB$ in
the form $\{x\} \cup Y.$

In order to prove that $\iat_\tau(N)$ is definable by
means of monadic second-order logic, we shall define
by means of monadic second-order logic the sets of
IA-automorphisms $\iat^+_\tau(N)$ and
$\iat^-_\tau(N).$ Here $\iat^+_\tau(N)$ denotes the
set of IA-automorphisms of the form
\begin{align*}
\alpha x &=x,\\
\alpha y &=y d_y, \quad y \in Y,
\end{align*}
where $d_y \in \str Y{}',$ and $\iat^-_\tau(N)$
the set of IA-automorphisms of the form
\begin{align*}
\alpha x &=x,\\
\alpha y &=y [x,a_y], \quad y \in Y,
\end{align*}
where $a_y \in \str Y.$

Clearly, $\iat_\tau(N)$ is a direct product
of $\iat^+_\tau(N)$ and $\iat^-_\tau(N).$

We shall call extremal involutions whose action on
$\cB$ is canonical {\it basis} extremal involutions; a
basis extremal involution taken an element $b$ from $\cB$
to the inverse will be denoted by $\f_b.$ The
automorphisms of $N$ which act on $\cB$ as
permutations will be called {\it basis permutations.}
By Lemma \ref{SymmsBasics} (b) $\f$ is a basis
extremal involution if and only if $\f$ is extremal
modulo $\ia N,$ commutes with $\theta$ and $\f B
\f \subseteq B^{\pm 1}.$ Similarly, the basis permutations are those
elements in the normalizer of $B$ which commute with
$\theta.$

{\it I. Characterization of $\iat^+_\tau.$}
The natural superset of $\iat^+_\tau$ is the set $C$
of IA-automorphisms in the centralizer of the basis
extremal involution $\f_x.$  By Lemma
\ref{PMs-of-an-Ext}, if an IA-automorphism $\gamma$
commutes with $\f_x,$ then $\gamma$ acts on $\cB
=\{x\} \cup Y$ as follows
\begin{alignat}2
&\gamma x &&= x[t,x],\\
&\gamma y &&=y d_y, \quad y \in Y, \nonumber
\end{alignat}
where $t$ is an element of $\str Y$
and $d_y \in \str Y{}'$ for each $y \in Y;$
it is easily seen that every basis extremal
involution normalizes $C.$

Thus, we have to choose those automorphisms $\gamma$
with (\theequation) that have $d(\gamma)=[t,x]$ equal to $1.$

Let $b$ be an element of $Y=\cB \setminus \{x\}$ and $\f_b$
the corresponding basis extremal involution.
For an IA-automorphism $\gamma_b=\sqrt{\mathstrut \f_b
\gamma\inv \f_b \gamma}$ obtained from an element $\gamma
\in C,$ the word $d(\gamma_b)$ has as a word in the
letters $\cB$ only occurrences of $x$ and $b.$ Really,
write the action of $\gamma$ on $\cB$ in the
following form
\begin{alignat*}2
&\gamma x &&= x [b^k u,x],\\
&\gamma y &&= y [b,e_y] f_y, \quad y \in Y \setminus \{b\},\\
&\gamma b &&= b [b,e_b] f_b,
\end{alignat*}
where $u \in \str{Y\setminus \{b\}},$ $e_y,e_b \in
\str{Y \setminus \{b\}}$ and $f_y,f_b$ are in
$\str{Y\setminus \{b\}}'.$ One easily
verifies that
\begin{alignat}2
&\f_b \gamma\inv \f_b \gamma x &&= x [b^k,x]^2,\\
&\f_b \gamma\inv \f_b \gamma y &&= y [b,e_y]^2, \quad y \in Y \setminus \{b\},
\nonumber\\
&\f_b \gamma\inv \f_b \gamma b &&= b  f_b^2. \nonumber
\end{alignat}

We can therefore conclude that {\it $\gamma \in C$
preserves $x$ if and only if for each $b \in Y$ so
does $\gamma_b.$}

Let us now fix $b.$ The automorphism $\gamma_b$ is an
element of the set
$$
D=D(b) = \{\delta \in C : \f_b \delta \f_b=\delta\inv  \}
$$
(according to Lemma \ref{PMs-of-an-Ext} $D$ is
equal to the set of square roots of the elements of
the form (\theequation)). Then in order to explain
when $\gamma_b$ preserves $x,$ it suffices to obtain a
general characterization of the subgroup of all elements
of $D$ which preserve $x.$

A `transvection' $U_b \in \aut N$ such that
\begin{align*}
U_b x &=xb,\\
U_b y &=y, \quad y \in Y
\end{align*}
can be modulo $\ia N$ characterized as one
of the automorphisms $U \in \aut N$ with
\begin{align*}
& U \tau U\inv = \tau \tau_b,\\
& U \nu U\inv =\nu, \quad \nu \in B \setminus \{\tau \}.
\end{align*}

That is all we need at the moment, since
we are going to act by $U_b$ on a set of IA-automorphisms
(namely on $D$) by conjugation:
$$
(U_b \beta) \alpha (U_b \beta)\inv = U_b \alpha U_b\inv
$$
for any $\alpha,\beta \in \ia N.$

We shall use below a family of automorphisms
$S=S(b)$ such that any member of $S$ fixes
all elements in $\cB \setminus \{b\}$
and takes $b$ to an element of the form
$b f,$ where $f \in \str{Y \setminus \{b\}}'.$

Let us obtain a description of the members of $S.$
Suppose $\delta$ is an automorphism from $D$ such that
\begin{align*}
& \delta x =x [b^k,x],\\
& \delta b = b f,
\end{align*}
where $f \in \str{Y \setminus \{b\}}'.$ One then
readily verifies that the automorphism $\delta^* = U_b
\delta\inv U_b\inv \delta$ takes $x$ to $xf$ and fixes
all elements of $Y.$ Therefore $\pi_{x,b} \delta^*
\pi_{x,b},$ where $\pi_{x,b}$ is the basis permutation
taking $x$ and $b$ to each other and preserving all
other elements of $\cB,$ is in $S.$

We claim now that the subgroup $D_0$ of the elements
of $D$ which stabilize $x$ can be described as a unique
subgroup $E$ of $D$ satisfying the following conditions (a-c):
\begin{itemize}

\item[(a)] $D$ is a direct product of $E$
and the subgroup generated by conjugation
by the basis element $b$: $D=E \times \str{\tau_b};$

\item[(b)] any basis permutation fixing both
$x$ and $b$ (or equivalently commuting with
both $\tau_x$ and $\tau_b$) normalizes $E;$

\item[(c)] $E \supseteq S.$
\end{itemize}

It is quite clear that $D_0$ satisfies the
conditions (a-c). Let us prove the
converse. Assume the conditions
(a-c) hold for a subgroup $E$ of $D.$

It can be deduced from (a) that
for any $\delta \in D_0$ there is
an integer $k \in \Z$ such that
$$
\tau_b^k \delta \in E.
$$

To apply (b) we will prove

\begin{lem}
For any automorphism $\delta \in D_0$ there exist $\mu
\in D_0,$ a basis permutation $\pi$ which fixes $x$ and
$b,$ and $\s \in S$ such that
$$
\delta = \mu^\pi \mu\inv \s
$$
{\em(}$\mu^\pi$ denotes $\pi \mu \pi\inv${\em)}.
\end{lem}

This will imply that $\delta \in E$ and
hence $D_0 \subseteq E;$ therefore
$D_0=E,$ since $D=D_0 \times \str{\tau_b}.$
Indeed, as we noted above there is $m \in \Z$ such
that $\tau_b^m \mu \in E.$ We then have
$$
(\tau_b^m \mu)^\pi (\tau_b^m \mu)\inv \in E \Rightarrow \mu^\pi \mu\inv \in E .
$$
Using (c), we obtain that $\delta=\mu^\pi \mu\inv \s$
is in $E.$

Let us prove the latter Lemma. Write the (infinite)
set $Y \setminus \{b\}$ in the form
$$
Y \setminus \{b\}=\{y_{i,k} : i \in I, k \in \Z\},
$$
where the cardinality of $I$ is equal to $\rank N.$

Suppose that
$$
\delta y_{i,k}=y_{i,k} [b,z_{i,k}],\quad i \in I,\quad k \in \Z,
$$
where $z_{i,k} \in \str{Y \setminus \{b\}},$
and consider a basis permutation $\bmpi$ which fixes
$x$ and $b,$ and which acts on $Y \setminus \{b\}$
as a permutation with infinite cycles:
$$
\bmpi y_{i,k}=y_{i,k+1}
$$
for every $i \in I$ and $k \in \Z.$

Therefore if $\mu \in D_0$ is defined as follows
\begin{alignat*}3
&\mu &&y_{i,k} &&=y_{i,k} [b,t_{i,k}], \\
&\mu &&b       &&=b,
\end{alignat*}
then the equation $\delta=\mu^{\bmpi} \mu\inv\sigma$ holds
for a suitable $\s \in S$ if, for instance,
$$
\bmpi t_{i,k-1} t_{i,k}\inv = z_{i,k} \quad \forall i \in I\quad \forall k \in \Z.
$$
An easy way to satisfy the latter equations is to set
all the elements of the form $t_{i,0}$ equal, say, to $1,$
and to define other elements $t_{i,k}$ (for any
given $i$) by induction: to define the elements
$d$ with positive second indices one can use the formula
$$
t_{i,k}=z_{i,k}\inv\bmpi t_{i,k-1}, \quad k \ge 1
$$
and for $d$ with negative second indices the formula
$$
t_{i,k-1}=\bmpi\inv (z_{i,k} t_{i,k}), \quad k \le 0.
$$

{\it II. Characterization of $\iat^-_\tau$.} As above
we start with a natural superset of $\iat^-_\tau$
related to the involution $\f_x$: the set of
IA-automorphisms
$$
L=\{ \lambda \in \ia N : \f_x \lambda \f_x =\lambda\inv \}.
$$
By Lemma \ref{PMs-of-an-Ext}
a member $\lambda$ of $L$ acts on the basis $\cB$
as follows:
\begin{align}
\lambda x &=xc,\\
\lambda y &=y [x,a_y], \quad y \in Y, \nonumber
\end{align}
where $c \in \str{Y}'$ and $a_y \in \str Y.$
This time the problem of determining whether $c=c(\lambda)$
in (\theequation) is equal to 1 is solved
quite easily. It turns out that $\iat^-_\tau$ is
a subset of the elements of $L$ such
that conjugation by one more involution $\psi$
take them to the inverse:
\begin{equation}
\iat^-_\tau=\{\lambda \in L : \psi \lambda \psi =\lambda\inv\}.
\end{equation}
This $\psi$ can be modulo $\ia N$ characterized as
a basis permutation with respect to
a basis $\{b,xb\} \cup (Y \setminus \{b\})$ of $N,$
where $b \in Y$:
\begin{itemize}
\item $\psi^2=\id;$
\item $\psi \tau_b \psi=\tau_{x} \tau_b$ (this
easily implies that $\psi \tau_x \psi =\tau_x\inv$);
\item $\psi \tau_y \psi =\tau_y$ for each
$y \in Y\setminus \{b\}.$
\end{itemize}

As usual without loss of generality we may
suppose that our $\psi$ acts on the basis
$\cB$ such that it moves $x$ to the
inverse, $b$ to $xb$ and preserves
all elements in $Y \setminus \{b\}.$
Let us now prove (\theequation).
The involution $\psi$ inverts any commutator of the form
$[x,a],$ where $a \in \str Y{}.$ This easily
implies that $\psi \lambda \psi =\lambda\inv$
for each $\lambda \in \iat^-_\tau.$
Conversely, choose $\lambda \in L$ and suppose that
$\psi \lambda \psi =\lambda\inv.$ Compare the action
of both automorphisms $\psi \lambda \psi$ and
$\lambda\inv$ on the element $b.$ We have
\begin{align*}
\psi \lambda \psi b=\psi \lambda(xb)=\psi(xc(\lambda)b[x,a_b])=
x\inv\,\psi c(\lambda)\, xb [x,a_b]\inv=b\,\psi c(\lambda)\, [x,a_b]\inv;
\end{align*}
the latter element must be equal to $\lambda\inv b,$
that is, to $b[x,a_b]\inv.$ Hence $c(\lambda)=1,$
$\lambda \in \iat^-_\tau,$ and the theorem is proved.
\end{proof}

\section{\it It fixes all conjugations and some symmetry} \label{Delta}

Suppose $\Delta$ is an automorphism of the group $\aut
N.$ Consider a basis set $B$ of conjugations. Since
the property of being a basis set of conjugations can
be expressed by means of monadic second-order logic,
$\Delta$ takes $B$ to other basis set of conjugations
$B'.$ On the other hand, there is an automorphism $\s
\in \aut N$ such that $\s B=B'.$ Thus, if we will
follow $\Delta$ by conjugation by $\s\inv$ we obtain
an automorphism of the group $\aut N$ which fixes all
members of the set $B,$ and therefore all conjugations
from the subgroup $\inn N.$

This automorphism, say, $\Delta_1,$ preserves surely
the normalizer of $B.$ Therefore if $\thetast$
is a symmetry attached to $B,$ then the image of $\thetast$
under $\Delta_1$ is a conjugate of $\thetast$ by an IA-automorphism
$\beta$ (Lemma \ref{Only-Conjs-of-a-Symm}). Follow then
$\Delta_1$ by conjugation by $\beta\inv$ we get
an automorphism of the group $\aut N$ which
fixes all conjugations and the chosen
symmetry $\thetast.$

Let us denote the latter automorphism again
by $\Delta.$ We are going then to prove
that $\Delta$ preserves all $\iat$-automorphism
(Proposition \ref{DeltaFixesIAs} below). This
will enable us to prove the main result of
the paper.

\begin{Th}
The automorphism group of an infinitely generated
free two-step nilpotent group is complete.
\end{Th}

{\it Proof} (assuming Proposition \ref{DeltaFixesIAs}).
We saw above that the group $\aut N$ is centreless
(Proposition \ref{Centreless}); we prove now that
$\Delta$ acts trivially on $\aut N.$

Let $\theta'$ be an arbitrary symmetry in $\aut N.$
Then the automorphism $\Delta$ must preserve
$\theta',$ because the product $\theta'\thetast$
is an $\iat$-automorphism, it must preserve by
Proposition \ref{DeltaFixesIAs}, and because it
must preserve $\thetast$ by the construction:
$$
\theta'\thetast=\Delta(\theta'\thetast)=\Delta(\theta')\Delta(\thetast)=
\Delta(\theta')\thetast \Rightarrow \Delta(\theta')=\theta'.
$$

\begin{Lem} \label{IfDeltaPreservesConjs}
If an automorphism of the group $\aut N$ preserves
all conjugations, then it preserves all elements
of $\aut N$ modulo the subgroup $\ia N.$
\end{Lem}

\begin{proof}
Let $\Gamma$ be an automorphism of $\aut N$
preserving all conjugations. Let further
$\sigma$ be an automorphism of $N$ and
$z$ an arbitrary element of $N.$ We then
have
$$
\tau_{\sigma z}=\Gamma(\tau_{\sigma z})=
\Gamma(\sigma \tau_z \sigma\inv) =
\Gamma(\sigma) \tau_z \Gamma(\sigma\inv)=
\tau_{\Gamma(\sigma)z}.
$$
\end{proof}

Take an arbitrary element $\sigma$ in $\aut N.$
The automorphism $\sigma \thetast\sigma\inv$
is a symmetry and we have that
\begin{equation}
\Delta(\sigma \thetast \sigma\inv)=\sigma \thetast\sigma\inv.
\end{equation}
On the other hand, by Lemma \ref{IfDeltaPreservesConjs}
$\Delta(\sigma) =\sigma \eta_\sigma,$ where $\eta_\sigma \in
\ia N$ and hence
$$
\Delta(\sigma \thetast \sigma\inv)=\sigma \eta_\sigma \thetast \eta_\sigma\inv \sigma\inv=
\sigma \thetast \eta_\sigma^{-2} \sigma\inv.
$$
This along with (\theequation) implies that
$\eta_\sigma=\id.$ The latter means that $\Delta$ stabilizes
all elements of the group $\aut N,$ and therefore any
automorphism of the group $\aut N$ is inner. $\qed$

\begin{Prop} \label{DeltaFixesIAs}
The automorphism $\Delta$ fixes all $\iat$-automorphisms.
\end{Prop}

\begin{proof}
First note that each subgroup $\iat_{\tau},$ where
$\tau$ is in $B,$ the basis set of conjugations
chosen above, is invariant under the automorphism
$\Delta,$ since $\iat_\tau(N)$ is definable in $\aut
N$ with the parameters $\tau,\thetast$ and $B$ by
means of monadic second-order logic (Theorem
\ref{Def-of-IAtau}) and $\Delta$ preserves each of the
parameters.

Let $\cB=\cB(B,\thetast)$ denote the basis of $N$ uniquely
determined by a pair $(B,\thetast)$ (see
Section \ref{conjs&symms}).
Take $\tau \in B$ and suppose that $\tau=\tau_b,$ where
$b \in \cB.$ Let $J_\tau$ denote the subgroup of
all $\iat$-automorphisms which fix each element
in $\cB \setminus \{b\}.$ Clearly any
$\iat$-automorphism $\eta$ is fully determined
by its projections on the subgroups
$J_\tau$ parallel to $\iat_\tau(N)$: for each $\tau \in B$ there
exists a unique $\eta_\tau \in J_\tau$ such
that
$$
\eta \eta_\tau\inv \in \iat_\tau(N).
$$
Therefore in order to complete the proof of the
Proposition it suffices to prove that $\Delta$
stabilizes the elements in all subgroups $J_\tau.$
These subgroups are all contained in the finitary
automorphism group $\operatorname{Aut}_{f,\cB}(N)$ of
$N$ consisting of all automorphisms of $N$ which fix
all but finitely many elements of $\cB.$

The automorphism group of a finitely generated free
two-step nilpotent group $N_k$ with a basis $\mathcal
X=\{x_1,x_2,\ldots,x_k\}$ can be generated by a set
consisting of one basis extremal involution, two basis
permutations (associated with $\mathcal X$) and an
automorphism $U$ of $N_k$ such that $U x_1 = x_1 x_2$ and $U
x_i = x_i,$ where $i \ge 2$ (see, e.g., \cite[pp. 272-273]{DFo}).
This implies that the group $\operatorname{Aut}_{f,\cB}(N)$
is generated by basis extremal involutions
(with respect to the basis $\cB$),
basis permutations acting on $\cB$ as finite
cycles and $U \in \aut N$ such that
\begin{align*}
U x &=x z,  \\
U b &=b,\quad \forall b \in \cB \setminus \{x\}
\end{align*}
where $x$ and $z$ are some distinct elements of $\cB.$

We claim that $\Delta$ preserves all elements
of the subgroup $\operatorname{Aut}_{f,\cB}(N).$
There are no problems with basis extremal involutions
and basis permutations, since they commute with $\thetast$:
if $\thetast \sigma \thetast=\sigma$ and
$\Delta(\sigma)=\sigma\eta,$ where $\eta \in \ia N$
(Lemma \ref{IfDeltaPreservesConjs}), then by
applying $\Delta$ to the both sides of the
first equation we have that $\eta=\id.$

It is easy to verify that the triplet
$(\tau_x \tau_z, \tau_x \thetast, B'),$ where
$B'=\{\tau_x\} \cup \{\tau_x \tau : \tau \in B, \tau \ne \tau_x\}$
codes (as it was described in Section \ref{conjs&symms})
the element $xz$ of $N.$ Let us denote by $\equiv$ the equivalence relation
defined by the conditions (PE1) and (PE2) (both
equivalent to group-theoretic ones by Theorem \ref{Def-of-IAtau})
from Section \ref{conjs&symms}. Then $U$ is the unique automorphism
$\sigma$ of $N$ satisfying the following conditions
with the parameters $\theta$ and $B$:
\begin{alignat*}2
& (\s \tau_x \s\inv, \s \theta \s\inv, \s B \s\inv) && \equiv
(\tau_x \tau_z, \tau_x \theta, B'),\\
& (\s \tau \s\inv, \s \theta \s\inv, \s B \s\inv) &&\equiv
(\tau, \theta,B),\quad \forall \tau \in B \setminus \{\tau_x\},
\end{alignat*}
and hence $\Delta(U)=U.$
\end{proof}

\end{document}